\documentclass[12pt,a4paper,oneside]{amsart}
\usepackage{amsfonts, amsmath, amssymb, amsthm, amscd, hyperref}
\usepackage[T2A]{fontenc}
\usepackage[cp1251]{inputenc}
\usepackage[english]{babel}
\usepackage[dvips]{graphicx}
\usepackage{subfig}
\usepackage{psfrag}
\usepackage{epstopdf}
\usepackage{cmap}
\usepackage{enumerate}
\usepackage[english]{babel}
\usepackage{amssymb}
\usepackage{amsmath}
\newtheorem{lemma}{Lemma}

\newtheorem{remark}{Remark}
\newtheorem{fact}{Theorem}
\newtheorem{corollary}{Corollary}
\newtheorem{conjecture}{Conjecture}

\theoremstyle{definition}
\newtheorem{definition}{Definition}
\newcommand{\la}{\lambda}

\newcommand{\e}{\varepsilon}

\def\frl{\forall}%
\newcommand{\men}{\leqslant}
\newcommand{\bol}{\geqslant}
\newcommand{\bra}{\langle}
\newcommand{\ket}{\rangle}
\newcommand{\B}{B}
\newcommand{\R}{\mathbb{R}}

\newcommand{\norm}[1]{\left\| #1 \right\|}

\def\vn{\mathop{\rm int}}

\def\Lin{\mathop{\rm Lin}}

\newcommand{\mglx}[1]{\rho_{X}\!\!\!\:\left( #1 \right)}\!%
\newcommand{\mco}[2]{\delta_{#1}\!\left( #2 \right)}%
\newcommand{\mcox}[1]{\delta_{X}\!\left( #1 \right)}%
\newcommand{\mgbx}[1]{\delta_{X}^{+}\!\left( #1 \right)}%
\newcommand{\lamx}[1]{\lambda^{-}_{X}\!\left( #1 \right)}%
\newcommand{\lapx}[1]{\lambda^{+}_{X}\!\left( #1 \right)}%
\newcommand{\mcoxi}[1]{\delta_{X}^{-1}\!\left( #1 \right)}%
\newcommand{\mcoi}[2]{\delta_{#1}^{-1}\!\left( #2 \right)}%
\newcommand{\prp}{\urcorner}%
\newcommand{\prooff}{\noindent {\bf Proof.}\\}%
\newcommand{\bbox}{\par\noindent\ensuremath{\Box}\par\noindent}%
\newenvironment{prf}
{\prooff}
{\bbox}

\def\B{\mathfrak{B}}%
\def\SS{\partial\B_1(o)}%
\def\SSS{\partial\B_1^*(o)}%
\def\BB{\B_1(o)}%
\newcommand{\reff}[1]{(\ref{#1})}%

\begin{document}

\title{Modulus of supporting convexity and supporting smoothness}

\author{G.M. Ivanov{$^{*}$}}

\address{Department of Higher Mathematics, Moscow Institute of Physics and Technology,  Institutskii pereulok 9, Dolgoprudny, Moscow
region, 141700, Russia}
\address{
National Research University Higher School of Economics,
School of Applied Mathematics and Information Science,
Bolshoi Trekhsvyatitelskiy~3, Moscow, 109028, Russia}
\thanks{{$^{*}$}Supported by the Russian Foundation for Basic Research, grant 13-01-00295.}
 \email{grimivanov@gmail.com}
\begin{abstract}
We introduce   the moduli of the supporting convexity and the supporting smoothness of a Banach space, 
which characterize the deviation of the unit sphere from an arbitrary supporting hyperplane.
We show that the modulus of supporting smoothness, the Bana{\'s} modulus, and the modulus of smoothness are all equivalent at zero, the modulus of supporting convexity is equivalent at zero to the modulus of convexity. We prove a Day--Nordlander type result for these moduli.    
\end{abstract}

\maketitle

\section{Introduction}
The properties of a Banach space are completely determined by its unit ball.
The geometry of the unit ball of a Banach space $X$ may be described, for instance, using the properties of some moduli attached to $X.$ (For example, the moduli of convexity, of smoothness, Milman's moduli, etc.)
The aim of this paper is to introduce and explore  some new type of moduli, 
which  characterize the deviation of the unit sphere from an arbitrary supporting hyperplane.

In the sequel we shall need some additional notation.
Let $X$ be a real Banach space. For a set $A \subset X$ by $  \partial A, \,\vn A$ 
we denote the boundary and the  interior of  $A.$
We use $\bra p,x \ket$ to denote the value of a functional $p \in X^*$ at  a vector $x \in X.$
For $R>0$ and $c \in X$ we denote by $\B_R(c)$ the closed ball with center $c$ and radius $R,$
by $\B_R^{*}(c)$ we denote the ball in the conjugate space.
By definition, put $J_1(x) = \{p \in \SSS :\, \bra p, x \ket\ = \norm{x}\}.$
For  convenience, the length of segment  $ab$ is denoted by  $\norm{ab},$
i.e., $\norm{ab} = \norm{a-b}.$ 

We say that $y$ is {\it quasiorthogonal}  to the vector $x \in X\setminus \{o\}$  
and write  $y\urcorner x$ 
if there exists a functional $p \in J_1(x)$  such that  $\bra p, y \ket = 0.$
Note that  the following conditions are equivalent: \\ \noindent
 -- $y$  is quasiorthogonal to  $x$ \\ \noindent 
 --  for any $\la \in \R$ the vector $x+ \la y$ 
 lies in the supporting hyperplane to the  ball  $\B_{\norm{x}}(o)$  at $x;$ \\ \noindent
 -- for any $\la \in \R$ the following inequality holds  $\norm{x + \la y} \bol \norm{x};$\\ \noindent
 --  $x$ is orthogonal to  $y$ in the sense of Birkhoff--James (\cite{DiestelEng}, Ch. 2, \S 1). 

Let
\begin{equation*}
\delta_X(\e) = \inf \left\{ 1 - \frac{\|x + y\|}{2}:\ x,y\in\B_{1}(0),\ \|x -y\| \ge \e\right\}
\end{equation*}
and 
\begin{equation*}
\mglx{\tau} = \sup \left\{\frac{\|x + y\|}{2} + \frac{\|x - y\|}{2} - 1 : \, \|x\| = 1, \|y\| = \tau \right\}.
\end{equation*}
The functions $\delta_X(\cdot): [0,2] \to [0, 1]$ and $\mglx{\cdot}: \R^+ \to \R^+$
 are  referred to as the moduli of convexity and smoothness of $X$ respectively.

Let $f$ and $g$  be two non-negative functions, each one defined on a segment $[0, \e].$
We shall consider  $f$ and $g$ as {\it equivalent at zero,} denoted by $f(t) \asymp g(t)$ as $t \to 0,$ 
if there exist positive constants 
$a,b,c,d,e$ such that $a f(bt) \men g(t) \men c f(dt)$ for $t \in [0,e].$

The rest of this paper is organized as follows.
In Section 2 we prove several simple technical lemmas,
in Section 3 we  introduce the definitions of the modulus of supporting convexity and the modulus of supporting smoothness and consider their basic properties, in Section 4 we show these modulus are equivalent at zero  to the modulus of convexity and smoothness respectively,  in Section 5 we prove that the moduli of  smoothness, of supporting smoothness and the  modulus  of Bana{\'s} are all equivalent at zero, and, finally, 
in Section 6 we prove some estimates for these moduli concerning the maximal value of the Lipschitz
constant for the metric projection operator onto  a hyperplane.

The author is grateful to professor G.E. Ivanov for constant attention to this work.

\section{Technical results}
In this section we prove several simple technical results.

The proof of the  next lemma is trivial.
\begin{lemma}\label{lemma o svyaznosti}
  Suppose the set $\B_1(o)\setminus \vn \B_r(o_1)$ is nonempty. Then it is arcwise connected.
\end{lemma}
\begin{lemma}\label{lemma_o_pleche}
Let $X_2$ be a two-dimensional Banach space.
Suppose $a,b, c , d \in \SS$ and the segments
 $ab,cd$ intersect in point  $x.$ 
Then the following inequality holds
$$
\min\{\norm{cx},\norm{xd}\} \men \max\{\norm{ax}, \norm{xb}\}.
$$
\end{lemma}
\begin{prf}
Assume the converse.
Then for some $\e > 0$ we get  $\min\{\norm{cx}, \norm{xd}\} > \max\{\norm{ax}, \norm{xb}\} + \e = r.$
Since the segment $ab$ belongs to  $\vn \B_{r}(x)$ and separates it into two parts, then
we cannot connect points  $c, d$  in  $\BB \setminus \vn \B_{r}(x).$ 
This contradicts Lemma  \ref{lemma o svyaznosti}. The lemma is proved.
\end{prf}
\begin{lemma}\label{UVO lemma ozen}
Let $x,y \in X,\; x\neq o,\; p \in \SSS$ such that $\bra p, x\ket = \norm{x}.$ Then
\begin{equation}
    \norm{x+y} \men \norm{x} + \bra p, y \ket + 2\norm{x} \mglx{{\frac{\norm{y}}{\norm{x}}}}.
\end{equation}
\end{lemma}
\begin{prf}
By definition of the modulus of smoothness, we get
$$
\frac{1}{2} \left( \frac{ \norm{x + y}}{\norm{x}} + \frac{\norm{x -y}}{\norm{x}} \right) - 1 \men \mglx{\frac{\norm{y}}{\norm{x}}}.
$$
Multiplying both sides by  $2\norm{x},$    after some transformations we obtain:
\begin{gather*}
\norm{x + y} \men 2\norm{x} - \norm{x-y} + 2 \norm{x} \mglx{\frac{\norm{y}}{\norm{x}}} \men\\
2\norm{x} + \bra p, y - x \ket + 2\norm{x} \mglx{\frac{\norm{y}}{\norm{x}}} = 
\norm{x} + \bra p, y \ket + 2\norm{x} \mglx{\frac{\norm{y}}{\norm{x}}}.
\end{gather*}
\end{prf}
\begin{lemma}\label{lemma_lapx_mgl_tehno1}
 For any vectors  $x,y\in X\setminus\{o\}$ the following inequality is true
$$
\norm{\frac{x}{\norm{x}} - \frac{y}{\norm{y}}}  \men  \frac{2\norm{x-y}}{\norm{x}}.
$$
\end{lemma}
\begin{prf}
Using the triangle inequality, we get
\begin{gather*}
	\norm{\frac{x}{\norm{x}} - \frac{y}{\norm{y}}} = 
	\norm{\left(\frac{x}{\norm{x}} - \frac{y}{\norm{x}}\right) + \left( \frac{y}{\norm{x}} - \frac{y}{\norm{y}}\right)} \men 
	\norm{\left(\frac{x}{\norm{x}} - \frac{y}{\norm{x}}\right)} + \norm{\left( \frac{y}{\norm{x}} - \frac{y}{\norm{y}}\right)} \men\\
	\men \frac{1}{\norm{x}} \norm{x-y}  + \norm{y} \left| \frac{1}{\norm{x}} - \frac{1}{\norm{y}}\right| \men \frac{2\norm{x-y}}{\norm{x}}.
\end{gather*}
\end{prf}

\section{Definitions and basic properties}
Let $x,y \in \SS$ be such that  $y\urcorner x.$ 
By definition, put  
$$\la_X(x,y,r) = \min{\{\la \in \R: \, \norm{x+ry - \la x} = 1\}}$$
for any  $r \in [0,1].$   
Denote  
\begin{equation*}
  \lamx{x,y,r} = \min\{\la_X(x,y,r), \; \la_X(x,-y,r)\};  
\quad
 \lapx{x,y,r} = \max\{\la_X(x,y,r), \; \la_X(x,-y,r)\}.
\end{equation*}

\begin{definition}\label{def_la1}
  For any $r \in [0,1]$ and $x \in \SS$ we define the {\it modulus of local supporting convexity} as
\begin{equation*}
   \lamx{x,r} = \inf {\lamx{x,y,t}},
\end{equation*}
and respectively, the  {\it modulus of local supporting smoothness} as
\begin{equation*}
 \lapx{x,r} = \sup{ \lapx{x,y,t}},
\end{equation*}
where we choose  $(y,t)$ such that 
$\norm{y} =1,$ $y \prp x,$ ${0 \men t \men r}$ 
to minimize (maximize) $\lamx{x,r}$  $(\lapx{x,r}).$
\end{definition}
It is clear that $\lamx{x,r} \men \lapx{x,r} \men 1.$
\begin{definition}\label{def_la2}
 For any $r \in [0,1]$  we define the {\it modulus of  supporting convexity} as
\begin{equation*}
   \lamx{r} = \inf {\lamx{x,t}},
\end{equation*}
and respectively, the  {\it modulus of  supporting smoothness} as
\begin{equation*}
 \lapx{r} = \sup{ \lapx{x,t}},
\end{equation*}
where we choose  $(x,t)$ such that 
$x \in \BB,$  ${0 \men t \men r}$ 
to minimize (maximize) $\lamx{r}$  $(\lapx{r}).$
\end{definition}

\begin{figure}[h]%
\center{
\psfrag{(1,0)}[1]{\raisebox{-1.8ex}{\hspace{1em}\tiny$(1,0)$}}
\psfrag{(0,0)}[1]{\raisebox{-1.4ex}{\hspace{1.3em}\tiny$(0,0)$}}
\psfrag{(-1,0)}[1]{\raisebox{-1.8ex}{\hspace{1em}\tiny$(-1,0)$}}
\psfrag{(0,1)}[1]{\tiny\hspace{1.4em}$(0,1)$}
\psfrag{la+}[1]{\tiny\hspace{4em}$\lapx{|a_1|}$}
\psfrag{la-}[1]{\tiny\raisebox{7.4ex}{\hspace{4.5em}$\lamx{|a_1|}$}}
\psfrag{f}[1]{\tiny\raisebox{2.4ex}{\hspace{-2em}$f(|a_1|)$}}
\psfrag{l}[1]{\tiny\hspace{1em}$\ell$}
\psfrag{S}[1]{\tiny\hspace{1em}$S$}
\includegraphics[height= 0.3\textwidth]{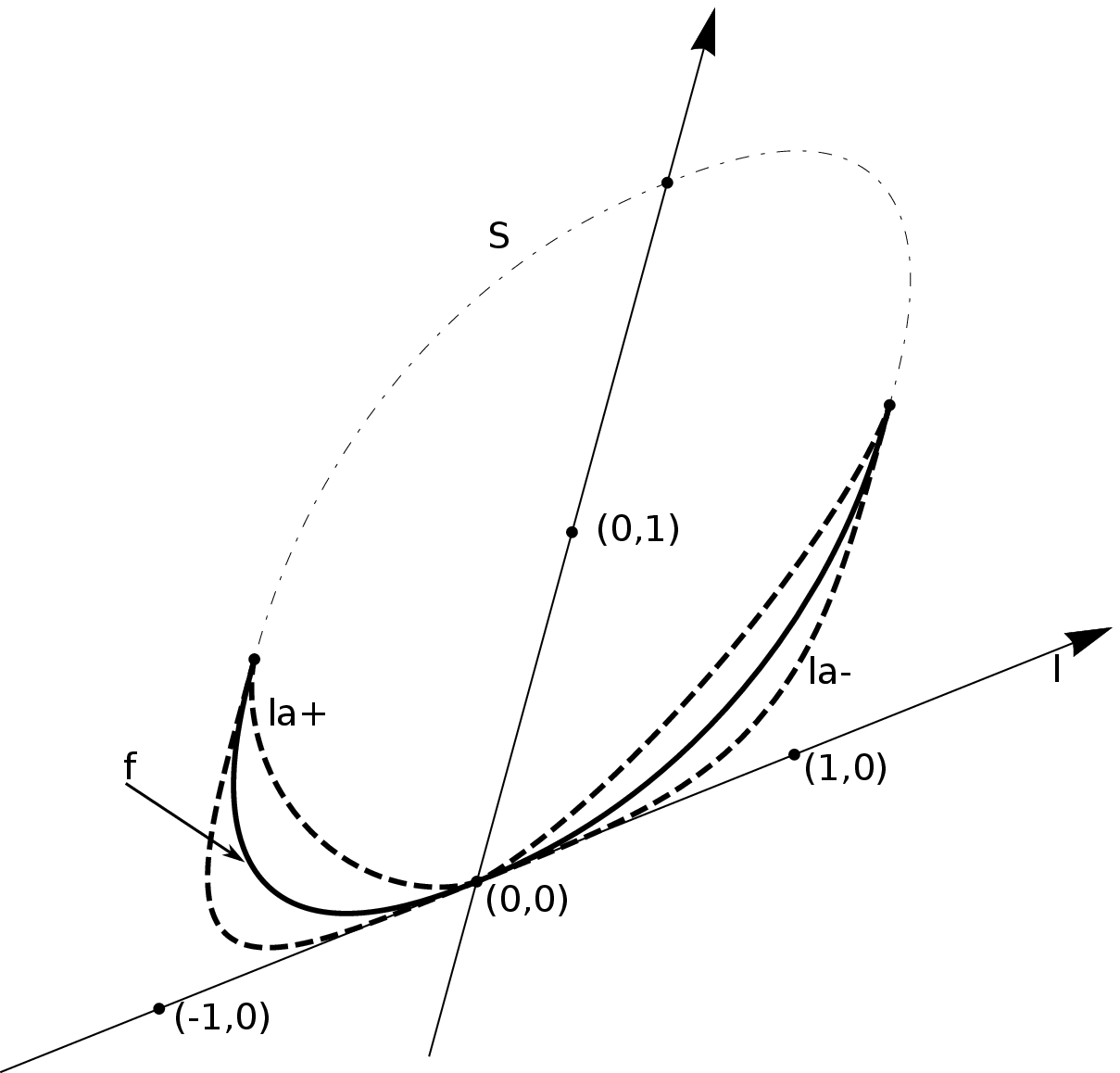}
}
\caption{Geometrical meaning of the $\lapx{r}, \lamx{r}.$ }
\label{figure_define_la}
\end{figure}

Let us explain the geometrical meaning of the moduli of supporting convexity and of supporting smoothness.
Fix $y, x \in \SS$ such that $y \prp x.$ 
Consider  the plane $L = \Lin\{y, x\}.$ We use $(a_1, a_2)$
to denote the  vector $a = a_1 y + a_2 x$ in this plane. 
The coordinate line $\ell=\{ (a_1,a_2)| a_1\in\R, a_2=0\}$ is a tangent to the unit "circle" $S=L\cap\partial\B_1(x)$. By the convexity of the ball, 
there is a convex function $f: [-1,1] \to \R$ such that 
for $a_1 \in [-1,1]$ the point $(a_1, f(a_1))$ belongs to the lower semicircle of $S$ (see Fig. \ref{figure_define_la}).
Hence for $a_1 \in [-1,1]$ the functions 
$\lamx{|a_1|}$ and $\lapx{|a_1|}$ are the lower and upper bounds to the $f(a_1)$ respectively, i.e.
the following inequalities hold $\lamx{|a_1|} \men f(a_1) \men \lapx{|a_1|}.$
\begin{lemma}\label{modop: osn_svoistva}
	Let $X$ be  an arbitrary Banach space, then:
\begin{enumerate}[(i)]
\item $\lapx{0}= \lamx{0} = 0;$
\item for any $r \in [0,1]$  the following inequality holds: $0 \men \lamx{r} \men \lapx{r}  \men r;$ 
\item for any  $0 < r_1 < r_2 < 1$ we have
\begin{eqnarray}
 \frac{r_2}{r_1} \lamx{r_1} \men \lamx{r_2} ,
\label{modop:1}\\
\lamx{r_2} - \lamx{r_1} \men \frac{r_2 - r_1}{1 - r_1};  
\label{modop:2}
\end{eqnarray}
\item the modulus of  supporting convexity is an increasing, continuous function on  $[0,1)$ and
moreover it is a strictly increasing function on the set  $\{r \in [0,1]: \lamx{r} > 0\};$
\item the modulus of  supporting smoothness is a strictly  increasing, convex and continuous function on  $[0,1]$
and furthermore $\lapx{1} = 1.$
\end{enumerate}
\end{lemma}
\begin{prf}
Let us introduce some notation.
Fix $x,y \in \SS$ such that $y \prp x,$ and real numbers $r_1, r_2$  such that 
$0 < r_1 < r_2 < 1.$ 
Let $z= x + y,$ $z_i = x + r_i y$ where $i=1,2.$ 
Let $y_1, y_2 \in \SS$ such that  
$y_iz_i \parallel ox$ and the intersection of the segment $y_iz_i$ and the ball $\BB$ 
is the point $y_i$ where $i=1,2.$
(see Fig. \ref{figure_lemma_141}).
\begin{figure}[h]%
\center{
\psfrag{x}[1]{$x$}
\psfrag{z}[1]{\raisebox{1ex}{\hspace{0.5em}$z$}}
\psfrag{o}[1]{$o$}
\psfrag{y}[1]{{\hspace{1em}$y$}}
\psfrag{y1}[1]{\hspace{1em}$y_1$}
\psfrag{y2}[1]{\hspace{1em}$y_2$}
\psfrag{z1}[1]{\hspace{1em}$z_1$}
\psfrag{z2}[1]{\hspace{1em}$z_2$}
\includegraphics[height= 0.3\textwidth]{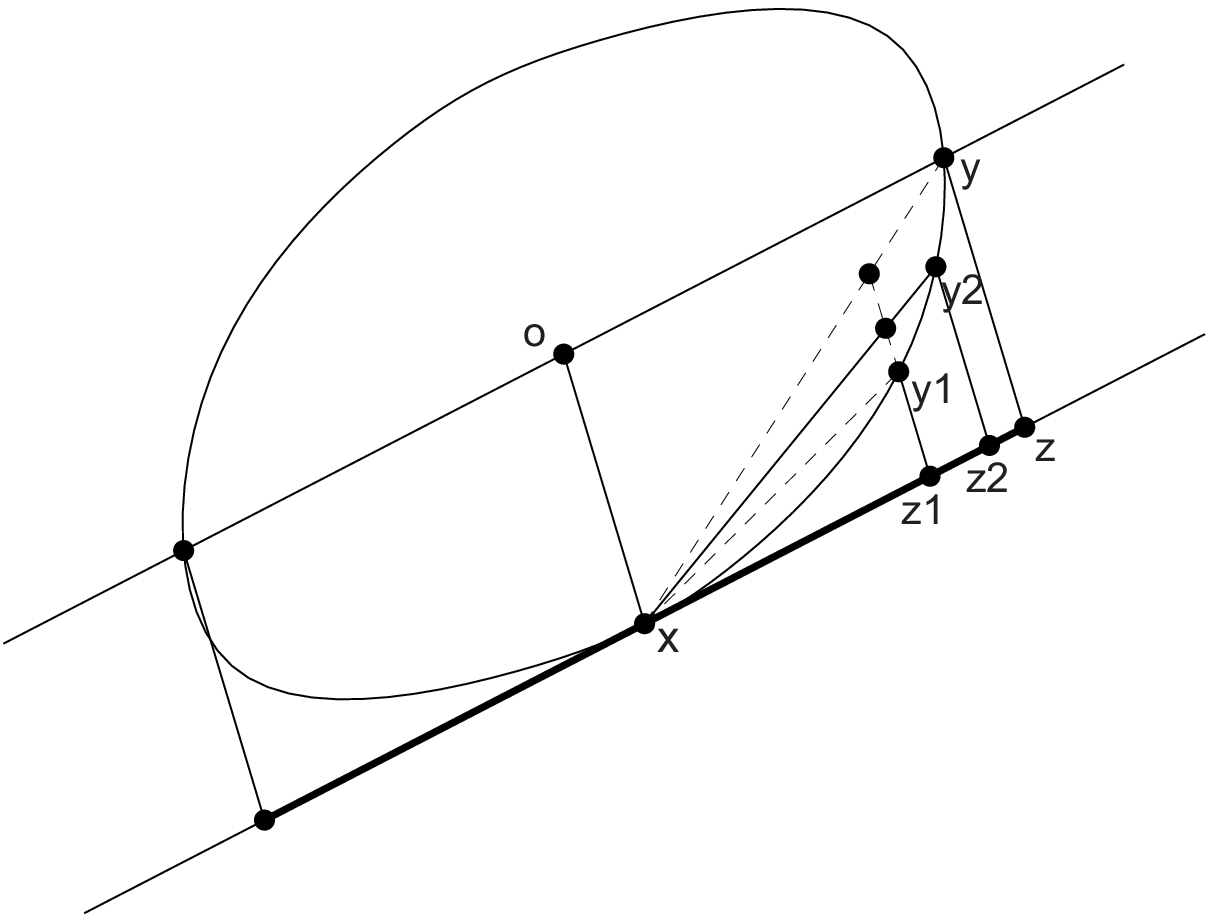}
}
\caption{ Illustration for Lemma \ref{modop: osn_svoistva}.}
\label{figure_lemma_141}
\end{figure} 
By  construction  $\norm{y_iz_i} = \la_X(x,y, r_i)$ where $i= 1,2.$
The reader will have no difficulty in showing that it is enough to prove all the assertions of this Lemma 
for  $\la_X(x,y,r).$
Now, let us prove the Lemma.\\ \noindent 
\vspace{-0.5cm}
\begin{enumerate}
	\item By the definitions, we have  $\lapx{0}= \lamx{0} = 0.$
	\item The first two inequalities of assertion (ii) are trivial.
	By similarity, we have $\la_X(x,y,r) \men r.$ Indeed, $y_1z_1\parallel zy$ and 
	$y_1z_1 \subset \triangle xyz.$ Taking the supremum we get assertion (ii).
	\item Taking into account that $\BB$ is convex, we get $y_1z_1 \subset xy_2z_2.$ 
	By construction we have that  $y_1z_1 \parallel z_2y_2.$   
	By the similarity, we get  $\norm{y_2z_2} \bol \frac{r_2}{r_1} \norm{y_1z_1},$
	 i.e.  $\frac{r_2}{r_1} \la_X(x,y,r_1) \men \la_X(x,y, r_2).$ 
	Taking the  infimum in $\la_X(x,y, r_2),$ we complete the proof of  inequality \reff{modop:1}.\\ \noindent
	By the convexity of the unit ball, we obtain that segment 	$y_2z_2$ lies in trapezoid 
	$y_1z_1zy.$ By construction  $y_2z_2 \parallel y_1z_1 \parallel yz.$ 
	By similarity, we get 
	$$\norm{y_2z_2} - \norm{y_1z_1} \men (1 - \norm{y_1z_1}) \frac{r_2 - r_1}{1 - r_1} \men \frac{r_2 - r_1}{1 - r_1}.$$
	Taking the infimum  in $\norm{y_1z_1} \to \lamx{r_1},$ we have 
	$\norm{y_2z_2} - \lamx{r_1} \men \frac{r_2 - r_1}{1 - r_1}.$
	This yields \reff{modop:2}.
	\item Assertion (iv) is the direct consequence of   assertion (iii).
	\item The function $\lapx{\cdot}$ is the supremum of the convex functions, therefore it's convex. 
	Since  $\lapx{\cdot}$ is a convex bounded function and 
	$\la_X(x,y,r)$ is continuous in $r,$ we obtain that $\lapx{\cdot}$ is continuous on $[0,1].$ 
	We will prove that  $\lapx{r} > 0$ on  $(0,1]$
	in Lemma \ref{modop:equiv_mgl_lapx2} below.   
	By this and  the equality $\lapx{0}=0$ and convexity of the modulus of supporting smoothness, 
	we get that it is a strictly increasing function. 
	The inequality $\lapx{r}\le r$ was proved in  assertion (ii).
	The equality $\lapx{1} = 1$ is the consequence of  inequality  \reff{la plus ozenka 1} at $r = 1,$
	which will be proved below. 
\end{enumerate}
\end{prf}

From Lemma \ref{modop: osn_svoistva} we have that in the definitions of the moduli of the supporting smoothness and supporting convexity one may choose $t = r.$
\begin{remark}\label{la_hilbert}
 Since any two plane central sections of the unit ball in a Hilbert space $H$ are equal, we have  
	\begin{equation*}
		\la^+_H(r)=\la^-_H(r) = \mco{H}{2r}  = 1 - \sqrt{1 - r^2}.
	\end{equation*}
\end{remark}
\section{Comparison of supporting moduli with the moduli of convexity and smoothness}
\begin{fact} \label{lemma o tolschine}
Let $X$ be an arbitrary Banach space. Then  $\lamx{\e} \asymp \mcox{\e}$ as $\e \to 0$ 
and for any  $r \in [0;1]$ :
  \begin{equation}
   \label{la ozenka snizu 6}
     \mcox{r} \men \lamx{r}  \men \mcox{2r}.
  \end{equation}
\end{fact}
\begin{prf}
1) By the definition of the modulus of supporting convexity for any  $\e > 0$
there exists a parallelogram  $xyzd$ such that 
$x, z \in \SS,$ the point $d$ lies in the segment  $xo$ and $\norm{xy} = r,$ $xy \prp ox,$ $\norm{yz}\men \lamx{r} + \e.$
 Therefore $\norm{od} = 1 - \norm{yz},$ consequently 
  $\mcox{r} =\mcox{\norm{zd}} \men \norm{yz} \men \lamx{r} + \e.$
Passing to the limit as $\e \to 0,$ we obtain the left-hand side of chain  \reff{la ozenka snizu 6}.\\
\noindent
2) Let us prove the right-hand side of  chain \reff{la ozenka snizu 6}.\\
Fix   $r \in (0,1)$ (if $r=0$ or $r=1$ the inequality is trivial).
 By the definition of  the modulus of supporting convexity for any  $\e > 0$
there exist points $a_\e, b_\e$ on the unit sphere such that 
$\norm{a_\e b_\e} \bol 2r$ 
and for the point  $c_\e = \frac{a_\e + b_\e}{2}$ the following inequality holds:
\begin{equation} \label{modop_la_sverh1}
1 - \norm{oc_\e} \men \mcox{2r} + \e.
\end{equation}
Let the ray $oc_\e$ intersect the unit sphere in a point $x.$ 
Denote by $l_1$ the supporting line to the unit sphere such that  
$l_1$ lies in the plane $oa_{\e}b_{\e}$ and $x \in l_1.$ 
Let $l_2$ be a line such that $l_1 \parallel l_2$ and $c_\e \in l_2.$  
Denote by  $f, g$ the points of intersections of $\SS$ and the line $l_2.$   
 From Lemma  \ref{lemma_o_pleche} it follows that  $\norm{f-c_\e} \bol r$ or $\norm{g-c_\e} \bol r.$ 
Without loss of generality, put
$\norm{g - c_\e} \bol r.$
Let $l_3$ be a line such that $l_3 \parallel oc_{\e}$ and $g \in l_\e.$  
By definition, we put  $y = l_3 \cap l_1 $ (see Fig. \ref{figure_lemma_142}).
\begin{figure}[h]%
\center{
\psfrag{x}[1]{$x$}
\psfrag{o}[1]{\hspace{0.21em}$o$}
\psfrag{f}[1]{\raisebox{-2ex}{$f$}}
\psfrag{y}[1]{\raisebox{-0.3ex}{$y$}}
\psfrag{g}[1]{\raisebox{0.7ex}{\hspace{0.6em}$g$}}
\psfrag{Ae}[1]{\raisebox{0.6ex}{\hspace{0.7em}$a_\e$}}
\psfrag{Be}[1]{\raisebox{1ex}{\hspace{1em}$b_\e$}}
\psfrag{Ce}[1]{\raisebox{1ex}{\hspace{1em}$c_\e$}}
\psfrag{l1}[1]{\raisebox{2ex}{$l_1$}}
\psfrag{l2}[1]{\raisebox{2ex}{$l_2$}}
\psfrag{l3}[1]{\hspace{0.7em}$l_3$}
\includegraphics[height= 0.3\textwidth]{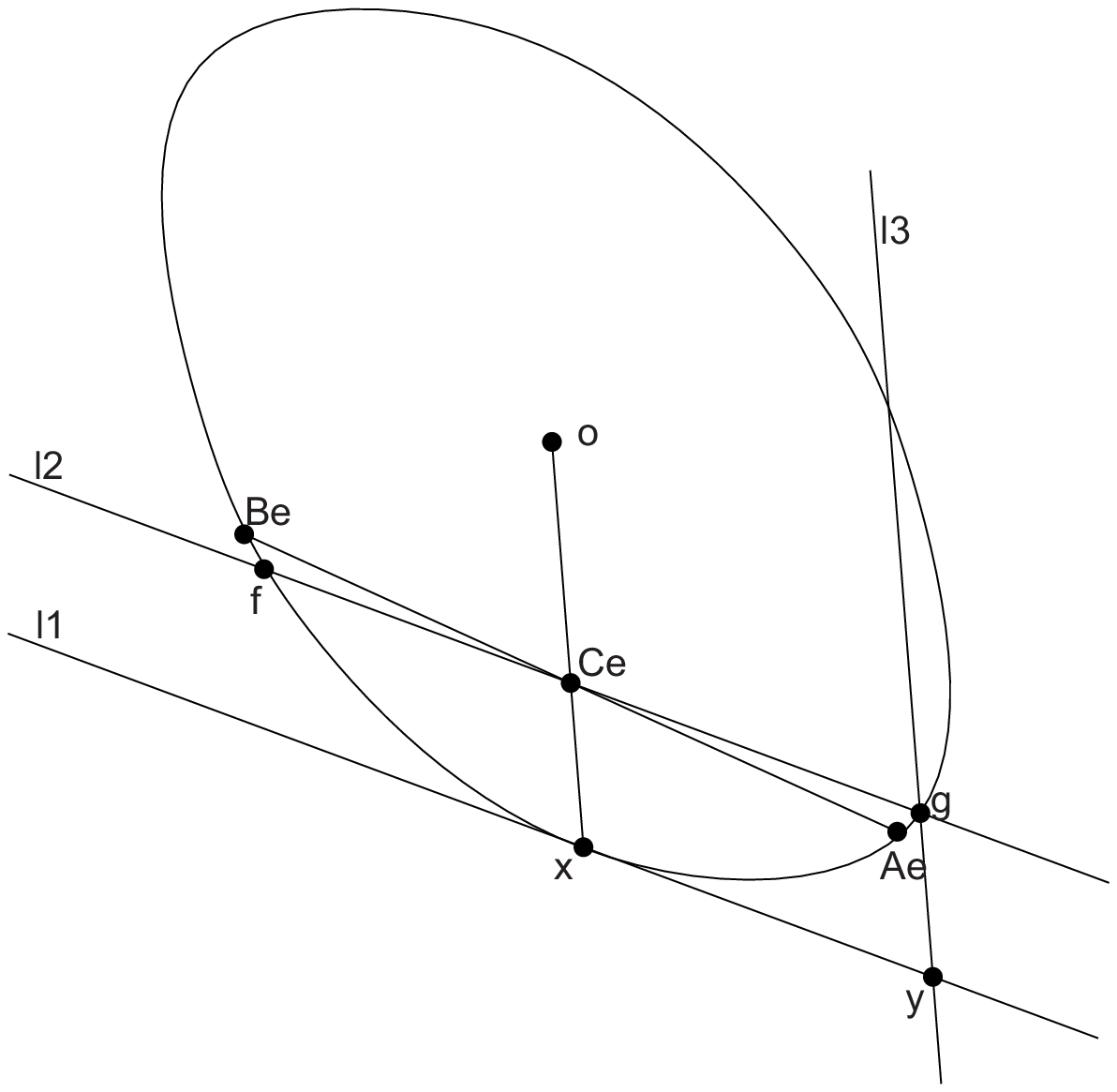}
}
\caption{Illustration for Lemma \ref{lemma o tolschine}.}
\label{figure_lemma_142}
\end{figure}
Then
\begin{equation*}
		\mcox{2r} + \e \bol \norm{c_{\e}x}  \bol \lamx{x, \frac{y-x}{\norm{y-x}}, \norm{y -x}} 
		\bol \lamx{x, \frac{y-x}{\norm{y-x}}, r} \bol \lamx{r},
\end{equation*}
i.e.,  $\mcox{2r} + \e  \bol \lamx{r}.$
Passing to the limit as  $\e \to 0,$ we complete the proof.
\end{prf}
\begin{lemma}\label{modop:equiv_mgl_lapx1}
    Let $r\in [0,\frac{1}{2}].$ Then
\begin{equation} \label{UVO la ozen}
  \lapx{r} \men  \mglx{2r}.
\end{equation}
\end{lemma}
\begin{prf}
 Denote $\la = \lapx{r}.$   
 Since $\lapx{r} \men r$ for any  $r \in [0,1],$ 
 then $\la \men \frac{1}{2}.$
 Let $\mu \in (0,\la).$ 
By  the Definitions \ref{def_la1}, \ref{def_la2}
there exist 
$x,y \in \SS$ such that $y\prp x$ and  $\lambda_X(x,y,r)=\mu'\in(\mu,\lambda),$
and consequently  $\norm{x+ry-\mu'x}=1$. 
Since $y\prp x$  there exists  $p\in J_1(x)=J_1(x-\mu'x)$ such that $\bra p,y\ket=0.$
 
Using Lemma \ref{UVO lemma ozen}, we get
\begin{gather*}
1 = \norm{x+ry-\mu'x} \men \norm{x-\mu'x}+\bra p,ry \ket + 2 (1 - \mu') 
\mglx{\frac{r}{1 - \mu'}} =  \\= 1 - \mu' + 2 (1 - \mu') \mglx{\frac{r}{1 - \mu'}}.
\end{gather*}
To complete the proof, it suffices to note that 
 $\mu' < \frac{1}{2},$ $\mglx{0} = 0$  and
 the modulus of smoothness is a convex function.
\end{prf}

\begin{lemma} \label{modop:equiv_mgl_lapx2}
	Let $r\in [0,1].$ Then
\begin{equation} \label{lemma_lapx_mgl_def}
 \mglx{\frac{r}{2}} \men \lapx{r}.
\end{equation}
\end{lemma}
\begin{prf}
Taking into account  the definition of the modulus of smoothness, it follows   that
for any  $\tau \in \left[0, \frac{1}{2}\right]$ and  $\e \in [0, \mglx{\tau})$ 
there exist $x$ and $y$ such that the following inequality is true
\begin{equation}\label{lapx_mgl_2}
\norm{x + \tau y} + \norm{x - \tau y} - 2 \bol 2(\mglx{\tau} - \e).
\end{equation}
Without loss of generality, we can assume that
 $\norm{x+ \tau y} \bol \norm{x-\tau y}$ (hence $\norm{x + \tau y} \bol 1$). 
Denote $u = \frac{x + \tau y}{\norm{x + \tau y}}, v = \frac{x - \tau y}{ \norm{x - \tau y}}.$
By Lemma \ref{lemma_lapx_mgl_tehno1}, we obtain 
\begin{equation} \label{lapx_mgl_1}
	\norm{u - v}  \men \frac{4\tau}{\norm{x+\tau y}};
\end{equation}
By the triangle inequality, we get
\begin{gather*}
\norm{u+v} \men \frac{2\norm{x}}{\norm{x+\tau y}} + 
\norm{x - \tau y} \left|\frac{1}{\norm{x +\tau y}} - \frac{1}{\norm{x -\tau y}}\right|\\ =
2 - \frac{1}{\norm{x +\tau y}}(\norm{x +\tau y} + \norm{x -\tau y} - 2).
\end{gather*}
Now, by inequality \reff{lapx_mgl_2}, we have that
\begin{equation}
\label{u+v}
\norm{u+v}  \men 2 - \frac{2(\mglx{\tau}-\e)}{\norm{x+\tau y}}.
\end{equation}

Let us consider the plane  $ouv.$ 
By  $\omega$ denote a point lying on the smallest arc $uv$ of the unit circle  such that 
the supporting line to the unit ball at $\omega$ is parallel to   $uv.$
Obviously,  either 
$\la_X\left(\omega, \frac{u-v}{\norm{u-v}}, \frac{\norm{u-v}}{2}\right) \bol 1 - \frac{\norm{u+v}}{2}$ 
or $\la_X\left(\omega, -\frac{u-v}{\norm{u-v}}, \frac{\norm{u-v}}{2}\right) \bol 1 - \frac{\norm{u+v}}{2},$
i.e. $\lapx{\frac{\norm{u-v}}{2}} \bol 1 - \frac{\norm{u+v}}{2}.$
Combining this with  inequalities  \reff{u+v}, we get
$$
  \frac{2(\mglx{\tau}-\e)}{\norm{x+\tau y}}  \men 2 \lapx{\frac{\norm{u-v}}{2}}.
$$
Now, by  inequality \reff{lapx_mgl_1}, we obtain 
$$
	\frac{2}{\norm{x+\tau y}} (\mglx{\tau} - \e) \men 2 \lapx{\frac{2 \tau}{\norm{x+\tau y}}} \men \frac{2}{\norm{x+\tau y}} \lapx{2\tau}.
$$ 
Multiplying both sides by  $\frac{\norm{x+\tau y}}{ 2}$ and passing to the limit as  $\e \to 0,$ 
we obtain \reff{lemma_lapx_mgl_def}.
\end{prf}

\begin{remark}
By Lemma  \ref{modop:equiv_mgl_lapx2} and the  properties of the modulus of smoothness, 
it follows that  $\lapx{r} > 0$ for all  $r>0.$
\end{remark}
By Lemmas \ref{modop:equiv_mgl_lapx1}, \ref{modop:equiv_mgl_lapx2} and the properties of the modulus of smoothness 
we have the following result.
\begin{fact} \label{modop:equiv_mgl_lapx3}
Let $X$ be an arbitrary Banach space. Then  $\lapx{\tau} \asymp \mglx{\tau}$ as $\tau \to 0$ 
and for any  $r \in [0, \frac{1}{2}]$:
$$
\mglx{\frac{r}{2}} \men \lapx{r} \men \mglx{2r}.
$$
\end{fact}

\section{Comparison with the Bana{\'s} modulus}
In the paper \cite{banas1}  J. Bana{\'s} defined and studied some new modulus of smoothness.
Namely,  he defined
$$\label{mod banass}
\mgbx{\e} = \sup \left\{ 1 - \frac{\norm{x + y}}{2}:\ x,y \in \BB,\ \norm{x -y} \men \e \right\}, \quad \e\in [0,2].
$$
The function $\mgbx{\cdot}$ is called the \textit{Bana{\'s} modulus}.
In the papers \cite{banas1, banas2, banas1990convexity, banas1997functions} 
several interesting results concerning this modulus
were obtained. 
Particulary, in   \cite{banas1}, J. Bana{\'s} proved that a space is uniformly smooth iff 
$
\frac{\mgbx{\e}}{\e} \to 0
$
 as  $\e \to 0.$ 
However, from the definition a space is uniformly smooth if and only if 
$
\frac{\mglx{\e}}{\e} \to 0
$
 as  $\e \to 0.$
This leads to the question: are the modulus of smoothness and the modulus of Bana{\'s} equivalent at zero?  
It is easy to check that there exist positive constant $a, b$ such that 
$\mgbx{t} \men a \mglx{bt},$ but the lower estimate of the modulus of Bana{\'s} in terms of the modulus of smoothness is unknown. 
In the next theorem  we prove that the modulus of Bana{\'s} and the modulus of supporting smoothness are equivalent at zero, so  Theorem \ref{modop:equiv_mgl_lapx3}  answers the  above question.
\begin{fact}\label{lemma_lapx_mgbx}
	Let  $X$ be an arbitrary Banach space. Then  $\mgbx{\e}\asymp\mcox{\e}$ as $\e \to 0$ and
	the following inequalities  hold:
  \begin{equation}
   \label{la plus ozenka 1}
    \mgbx{2r} \leq \lapx{r} \qquad \frl r \in \left[0, 1 \right]; 
  \end{equation}
	  \begin{equation}
   \label{la plus ozenka 2}
     \lapx{r} \leq 2\mgbx{3r} \qquad \frl r \in \left[0, \frac{2}{3}\right].
  \end{equation}
\end{fact}
\begin{prf}
1) First we shall prove  inequality \reff{la plus ozenka 1} for $r \in [0, 1).$\\
Let $a,b$ be  points of the  unit sphere such that $\norm{a-b} \men 2r.$
By $X_2$ denote the plane  $aob.$

There exists a point  $y_2$ of the unit sphere of the plane $X_2$   such that  
the supporting line $l_2$ to the unit ball  at this point is parallel to  $ab.$
By definition, put $y_1 = oy_2 \cap ab.$ 
There exists a point
$a_2$  in the  projection of the  point $a$ on $l_2$ such that  the segments $y_1y_2, aa_2$  
are equal in length and parallel.
The point $b_2$ is defined in  the same way, such that  $y_1y_2$ and $bb_2$ are parallel  
  (see Fig. \ref{figure_lemma_146}).
  Without loss of generality we assume that
 $\norm{y_2a_2} \men r < 1.$ \ 
Since the modulus of supporting smoothness is an increasing function, we have
 $\norm{y_1y_2} = \norm{aa_2} \men \lapx{y_2, \norm{y_2a_2}} \men  \lapx{y_2, r}.$ 
Taking the supremum, we obtain  inequality  \reff{la plus ozenka 1}.

Taking into account that the modulus of Bana\'s is a continuous and increasing function, we obtain  inequality 
\reff{la plus ozenka 1} for  $r=1.$

\noindent
2) Let us prove  inequality \reff{la plus ozenka 2}.\\
By the  definition of modulus of supporting smoothness  for any 
 $\e \in (0, \lapx{r})$ there exist\\ 
-- a point  $x\in \SS;$\\
--  a line $\ell_1$ supporting to the unit ball at point  $x;$\\
-- a point $y$ on $\ell_1$ and a point $z \in \SS$ such that\\
	$\norm{xy} = r,$ $\norm{yz} > 0,$ $zy \parallel ox$ 
	and
 $\la^+\!\!\left(x,\frac{xy}{\norm{xy}}, r\right)  = \norm{yz} > \lapx{r} - \e> 0.$ 

Let $\ell_2$ be a line parallel to $\ell_1$ such that $z \in \ell_2.$
Let $z, z_1$ be points of the intersections of line $\ell_2$ and $\SS.$
By $y_1$ denote  the projections of  $z_1$ on $\ell_1$ such that $z_1y_1 \parallel ox$ 
(see Fig. \ref{figure_Lemma147_2}).
\noindent
\begin{figure}[ht]%
\center{
\psfrag{x}[1]{$x$}
\psfrag{o}[1]{\raisebox{0.7ex}{$o$}}
\psfrag{y}[1]{$y$}
\psfrag{z}[1]{\raisebox{0.2ex}{\hspace{0.5em}$z$}}
\psfrag{y1}[1]{\raisebox{-3.5ex}{$y_1$}}
\psfrag{z1}[1]{\raisebox{1.5ex}{\hspace{-1em}$z_1$}}
\psfrag{d}[1]{\hspace{0.2em}$d$}
\psfrag{f}[1]{$f$}
\psfrag{e}[1]{\raisebox{0.3ex}{\hspace{0.5em}$e$}}
\psfrag{l1}[1]{\raisebox{2ex}{$\ell_1$}}
\psfrag{l2}[1]{\raisebox{2ex}{$\ell_2$}}
\includegraphics[scale=0.5]{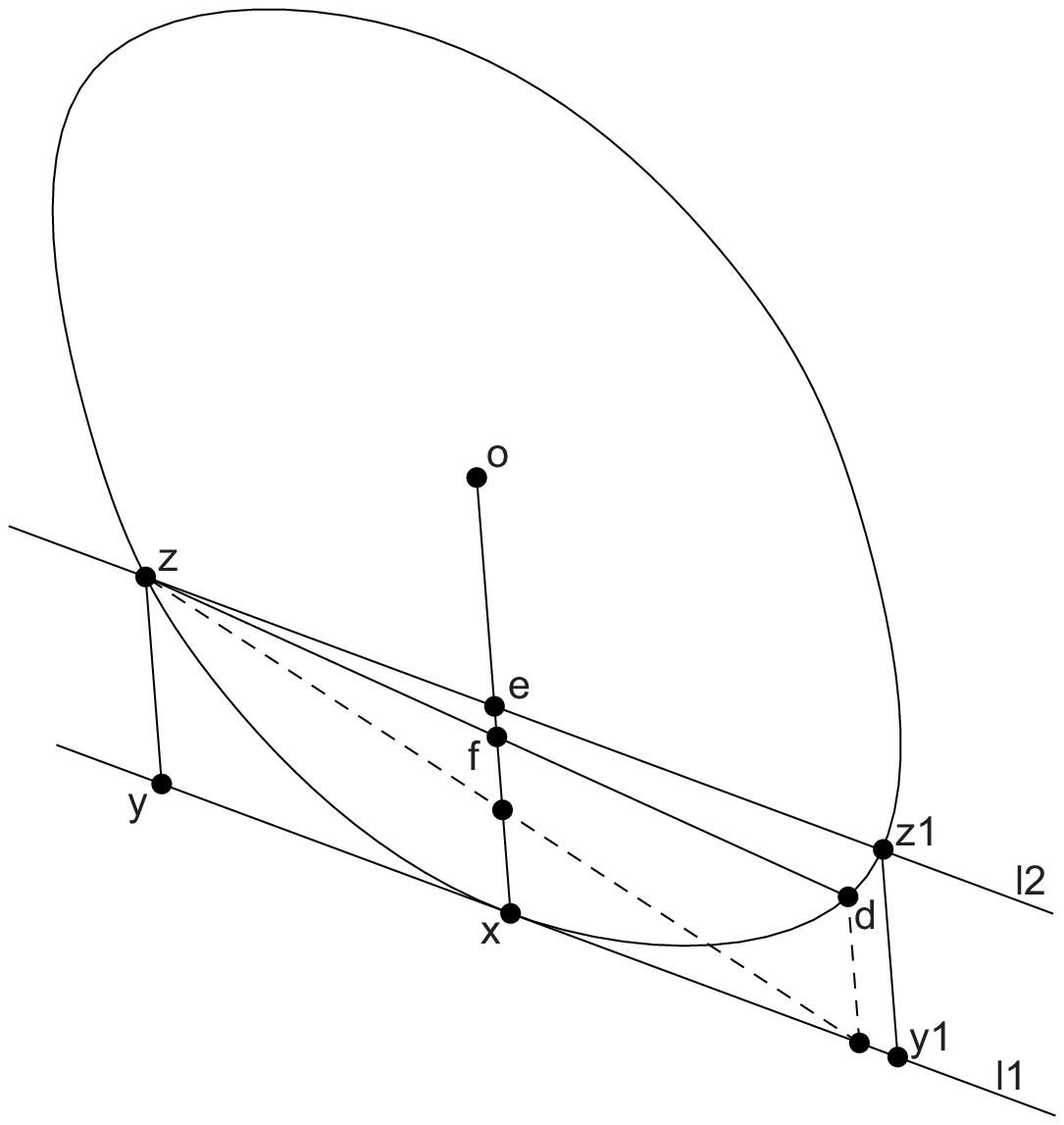}
}
	\caption{Illustration for the second part of Theorem \ref{lemma_lapx_mgbx}.}
		\label{figure_Lemma147_2}
\end{figure}

We shall prove that $\norm{zz_1} \bol 2r.$
In the converse case, $\norm{xy_1} < r.$ 
Note that if we fix $x,y \in \SS$ such that $y \prp x,$ 
then the function $\la^+\!\!\left(x,y, \cdot\right)$ is strictly increasing on the set of its positive values. 
Since $xy_1$ and $xy$  lie  on the same line and by to the definition of $\la^+,$ we obtain
$$\la^+\!\!\left(x,\frac{xy}{\norm{xy}}, r\right)\! = \norm{yz} = \norm{y_1z_1} \men 
\la^+\!\!\left(x,\frac{xy_1}{\norm{xy_1}}, \norm{xy_1}\right)\! < \!
\la^+\!\!\left(x,\frac{xy_1}{\norm{xy_1}}, r\right)\! = \! 
\la^+\!\!\left(x,\frac{xy}{\norm{xy}}, r\right)\!.$$ 
Contradiction. 
Consequently $\norm{zz_1} \bol 2r.$

By definition, put $e = ox \cap zz_1.$ 
By the continuity reasons there exists a point $d$ on the arc $z_1x$ of the unit sphere such that 
for the point $f = zd \cap ox $ the following equality holds $\norm{d-f} = \norm{f - z}.$  
Since  $\ell_1$ is a supporting line to the unit sphere, we have $\norm{xf} \bol \frac{\norm{yz}}{2}.$
Note that  $\norm{dz} \men 2(\norm{ze} + \norm{ef}) \men 3r.$ 
Combining the last two  inequalities, we get
 $$
	\mgbx{3r} \bol \norm{xf} \bol \frac{\lapx{r}-\e}{2}.
$$
Passing to the limit as  $\e \to  0,$ we obtain  inequality \reff{la plus ozenka 2}.
\end{prf}

From Theorems \ref{modop:equiv_mgl_lapx3} and \ref{lemma_lapx_mgbx} we have the following corollary.
\begin{corollary}
	Let  $X$ be an arbitrary Banach space, then  $\mgbx{\e}\asymp\mglx{\e}$ as $\e \to 0$ and
	the following inequalities  hold:
$$
	\frac{1}{2} \mglx{\frac{r}{6}} \men \mgbx{r} \men \mglx{r}, \quad r \in 
	\left[0, \frac{1}{2}\right].
$$
\end{corollary}

The  Day-Nordlander theorem (see \cite{Day-Nord}) asserts that $\mcox{\e} \men \delta_H(\e)$ for $\e \in [0,2],$
where $H$ denotes an arbitrary Hilbert space. 
On the other hand, repeating the arguments from the paper \cite{Day-Nord} we can show
that for any Banach space the following estimate is true $\delta^{+}_H(\e) \men \mgbx{r}$ for $\e \in [0,2].$
From this and  Theorems \ref{modop:equiv_mgl_lapx3}, \ref{lemma_lapx_mgbx}  
we obtain a Day--Nordlander type result for the moduli of supporting convexity and supporting smoothness:
\begin{corollary}\label{corol_daynord_lapx}
Let  $X$ be an arbitrary Banach space. Then
\begin{equation*}
\lamx{r} \men \la_{H}^{-}(r) = 1 - \sqrt{1 - r^2}= \la_{H}^{+}(r) \men \lapx{r} \qquad \frl r \in [0,1].
\end{equation*}

If at least one of these  inequalities turns into equality, then $X$ is  a Hilbert space.
\end{corollary}

\section{Estimates for Lipschitz constant for the metric projection onto a hyperplane}
Let us introduce the following characteristic of a space:
$$
    \xi_X = \sup\limits_{{\norm{x}=1,} \atop {\norm{y}=1}}\sup\limits_{p \in J_1(y)}{\norm{x - \bra p, x \ket y}}.
$$

Note that if $y \in \partial \B_1(0), \;p \in J_1(y),$ 
then the vector $(x - \bra p,x\ket y)$ is a metric projection of $x$ onto the hyperplane $H_p = \{x \in X: \bra p, x \ket = 0\}.$ 
So, $\xi_X = \sup_{y \in \B_1(o)}\sup_{p \in J_1(x)}{\xi_X^p},$ 
where $\xi_X^p$ is  half of diameter of a unit ball's projection onto the hyperplane $H_p.$ 
Therefore, $\xi_X$ is the maximal value of the Lipschitz constant for the metric projection operator onto a hyperplane.
Obviously, $\xi_X\le 2$ and $\xi_H = 1$ for a Hilbert space $H.$ 
\begin{fact}\label{th_chi_la}
For any Banach space $X$ the following inequality is true:
\begin{equation}
  \label{chi_E ozenka1}
 \frac{1}{1- \lamx{\frac{1 -\lamx{1}}{2}}} \men \xi_X \men \frac{1}{1 - \lapx{\frac{1 -\lamx{1}}{2}}}.
\end{equation}
\end{fact}
\begin{prf}
First let us introduce some notation.
Let  $x_0$ be an arbitrary point on the unit sphere. 
Let $l$ be a supporting line to the unit ball  at the point $x_0.$
Define $l_2$ as the line such that the following conditions  hold: \\ \noindent
a) $l_2 \parallel ox_0;$ \\ \noindent
b) $l_2 \cap l \neq \emptyset,$ by definition, put  $x_2 = l_2 \cap l;$ \\ \noindent
c)  $l$ a is supporting line to the unit ball at some point $y_2;$ \\ \noindent
d)  $\norm{y_2 x_2} \men 1.$ \\  \noindent
Let  $x_1$ be a point on segment $x_0 x_2$ such that $\norm{x_0 x_1} = 1,$ 
let $l_1$ be a line such that $x_1 \in l_1$ and $l_1 \parallel ox_0.$ 
By definition, put $y_1$ as the intersection point  of  line $l_1$ and the segment  $oy_2.$  
Let $b$ be a point on $\SS$ such that the segment $ob$  is parallel to  $x_0x_1.$ 
By construction, we have that  $x_0x_1bo$ is a parallelogram, 
therefore $b \in l_1$ and $y_1 \in x_1b.$
Let $a$ be the intersection point of the line $l_1 $
and the unit sphere such that  $a \in x_1y_1.$

From the intercept theorem, we have
 $\frac{\norm{x_0 x_2}}{\norm{o y_2}}=\frac{\norm{x_0 x_1}}{\norm{o y_1}}.$ 
 Therefore
\begin{equation} \label{gen ozenka}
\norm{x_0x_2} = \frac{1}{\norm{o y_1}} =\frac{1}{ 1 -\norm{y_1y_2}}.
\end{equation}
Since  $x_0x_1bo$ is a parallelogram, we get  $\norm{x_1b} =\norm{ox_0}= 1.$
By construction we have that $\norm{x_0x_1} =1.$ Therefore,  
\begin{equation} \label{chi_from_la}
\norm{ab} \men  1 - \lamx{1}.
\end{equation} 
Define $a_2$ as the projection of the point $a$ on $l_2$  such that $aa_2 \parallel oy_2.$
In the same way we define the point $b_2.$
Then the segments $y_1y_2, aa_2$ и $bb_2$ are parallel and equal in length 
(as parallel segments between two parallel lines). 
By the definition of the modulus of supporting convexity and by  inequality \reff{chi_from_la},
we obtain 
\begin{equation}\label{chi_from_la3}
\norm{y_1y_2} \men \lapx{\min\{\norm{a_2y_2}, \norm{y_2b_2}\}} \men \lapx{\frac{\norm{ab}}{2}} \men \lapx{\frac{1 -\lamx{1}}{2}}.
\end{equation}
Combining this and   equality \reff{gen ozenka}, 
we finally prove the right-hand side of inequality \reff{chi_E ozenka1}.

Let $\e$ be an arbitrary positive real number.
Note that we could choose a point  $x_0$ such that $\norm{x_1a} \men \lamx{1}+ \e,$
i.e. $\norm{ab} \bol 1- \lamx{1} - \e.$ 
Like in  \reff{chi_from_la3}, we obtain   
$$
	\norm{y_1y_2} \bol \lamx{\max\{\norm{a_2y_2}, \norm{y_2b_2}\}} \bol \lamx{\frac{\norm{ab}}{2}} 
	\bol \lamx{\frac{1 -\lamx{1} - \e}{2}}.
$$
Passing to limit as  $\e \to 0$ and using inequality \reff{gen ozenka}, 
we prove the left-hand side of  inequality  \reff{chi_E ozenka1}.
\end{prf}
\begin{remark}
  The estimate \reff{chi_E ozenka1}  is reached in case of a Hilbert space.
  The right-hand side of inequality  \reff{chi_E ozenka1} does not exceed 2,
  i.e. this estimate is not trivial.
  \end{remark}
\begin{conjecture}
 The right-hand side of inequality  \reff{chi_E ozenka1} becomes an equality  in case of  $L_p, \; p \in (1; +\infty)$. 
\end{conjecture}

\begin{figure}[ht]%
\center{
\psfrag{x0}[1]{\raisebox{1.7ex}{\hspace{0.8em}$x_0$}}
\psfrag{o}[1]{\hspace{0.6em}$o$}
\psfrag{x1}[1]{\raisebox{1.7ex}{\hspace{1em}$x_1$}}
\psfrag{x2}[1]{\raisebox{1.7ex}{\hspace{1.1em}$x_2$}}
\psfrag{y1}[1]{\raisebox{6.1ex}{\hspace{0.4em}$y_1$}}
\psfrag{y2}[1]{\hspace{1.4em}$y_2$}
\psfrag{a}[1]{\raisebox{5ex}{\hspace{1.6em}$a$}}
\psfrag{b}[1]{\raisebox{0.7ex}{\hspace{0.4em}$b$}}
\psfrag{a2}[1]{\hspace{1.3em}$a_2$}
\psfrag{b2}[1]{\hspace{1.4em}$b_2$}
\psfrag{l1}[1]{\raisebox{3.9ex}{\hspace{0.6em}$l_1$}}
\psfrag{l2}[1]{\raisebox{3.7ex}{\hspace{0.8em}$l_2$}}
\psfrag{l}[1]{\raisebox{1.7ex}{$l$}}
	\includegraphics{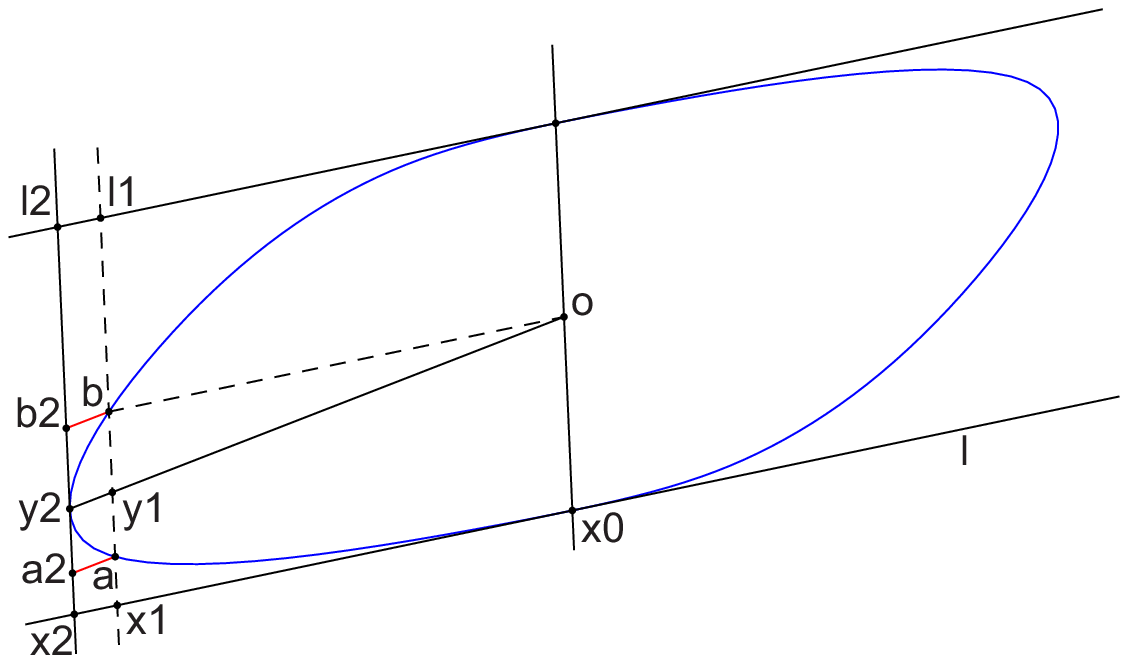}
}
	\caption{Illustration for Theorem \ref{th_chi_la}.}
	\label{fig_th_141}
\end{figure}

In the following lemma we obtain a lower estimate of the modulus of  supporting smoothness   
by the inverse function to the modulus of convexity. 
\begin{lemma} \label{lemma 1}
 For any  $r \in [0,1]$ the following inequalities  hold:
  \begin{equation} \label{la ozenka snizu}
   1 - \frac{1}{2} \mcoxi{1 - \frac{r}{2}}
    \men 
	1 - \frac{1}{2} \mcoxi{1 - \frac{r}{\xi_X}}
	 \men \lapx{r}.
  \end{equation}
\end{lemma}
\begin{prf}
The left-hand side of  inequality \reff{la ozenka snizu} is a straightforward consequence of the 
inequality  $\xi_X\le 2.$ 
Let us prove the right-hand side of inequality \reff{la ozenka snizu}.
In  case of $r=0$ it    is trivial.
Let $x_0$ be an arbitrary point on the unit sphere.
Define $H_x$ as a supporting hyperplane to the unit ball at the point $x_0.$  
Let $x_1$ be a point of the supporting hyperplane $H_x$ such that  
$\norm{x_0x_1} = r.$   
Denote the ray $\{ox_0 + \alpha x_0x_1: \; \alpha \bol 0 \}$  as $\ell.$ 
Let $l_1, l_2$ be the lines  parallel to  $ox_0$ such that \\ \noindent
a) $l_2$ is a supporting line to the unit ball at the point $y_2$ and 
 $l_2 \cap \ell = x_2;$ \\ \noindent
b) $l_1$ intersects the ray $\ell$ at    $x_1$ and 
intersects the unit sphere at points  $a, b.$ \\ \noindent
Let $y_1 = oy_2 \cap ab$ (see Fig. \ref{figure_lemma_146}).
\begin{figure}[ht]%
\center{
\psfrag{x0}[1]{\hspace{1.2em}$x_0$}
\psfrag{o}[1]{$o$}
\psfrag{x1}[1]{\hspace{1.2em}$x_1$}
\psfrag{x2}[1]{\hspace{1.2em}$x_2$}
\psfrag{y1}[1]{\hspace{1em}$y_1$}
\psfrag{y2}[1]{\hspace{1em}$y_2$}
\psfrag{a}[1]{\hspace{0.6em}$a$}
\psfrag{b}[1]{\hspace{-1em}$b$}
\psfrag{a1}[1]{\hspace{1em}$a_2$}
\psfrag{b1}[1]{\hspace{0.8em}$b_2$}
\psfrag{l1}[1]{\raisebox{1.7ex}{\hspace{-1em}$l_1$}}
\psfrag{l2}[1]{\raisebox{1.7ex}{\hspace{-1em}$l_2$}}
\psfrag{ll}[1]{$\ell$}
\includegraphics[height= 0.3\textwidth]{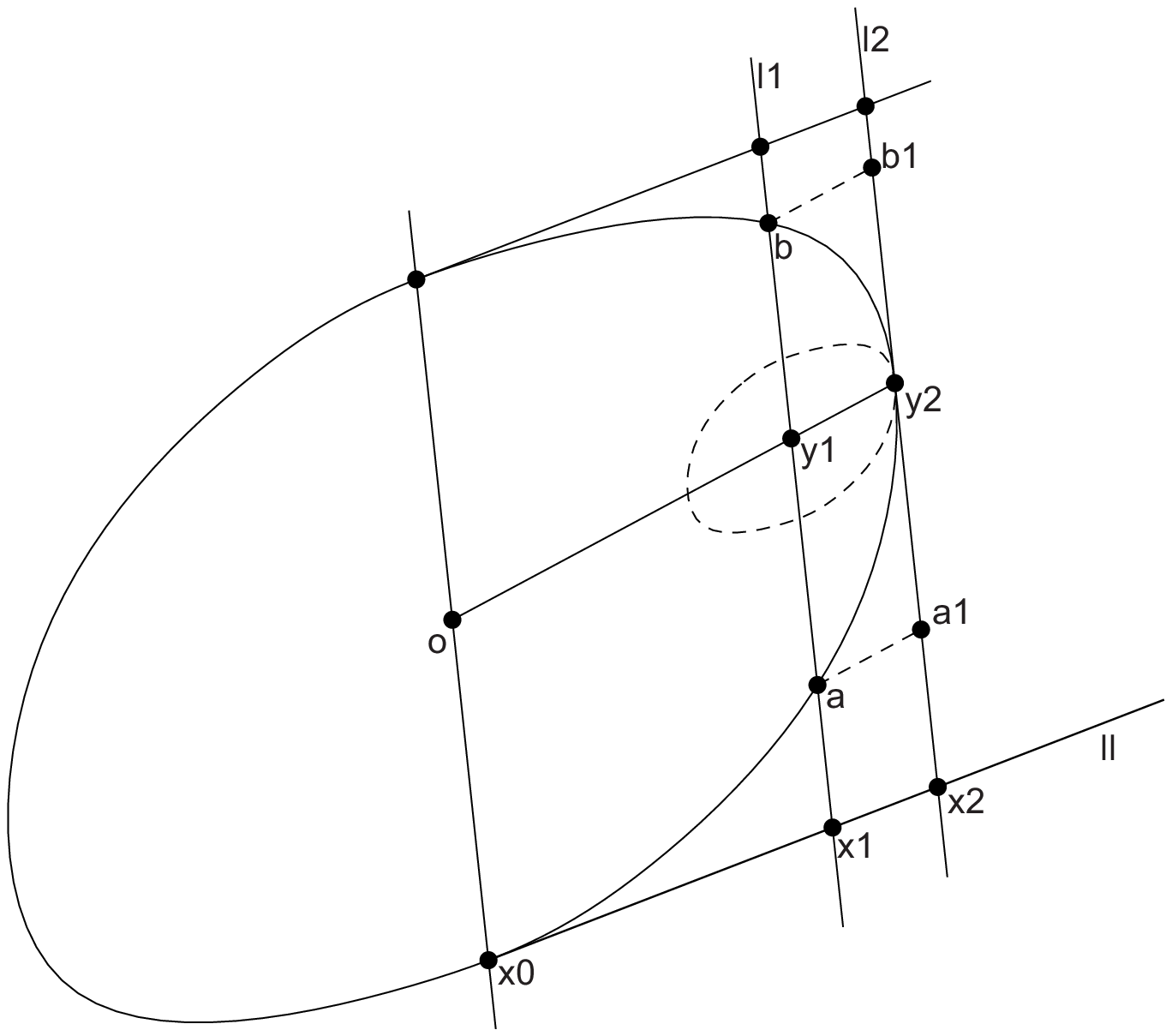}
}
\caption{Illustration for Lemma \ref{lemma 1} and for the first part of Theorem \ref{lemma_lapx_mgbx}.}
\label{figure_lemma_146}
\end{figure}
 By the definition of $\lapx{r}$ and 	since the unit ball is centrally symmetric,
 we get 
$\norm{ab} \bol 2(1- \lapx{r}).$ 
Obviously,  $\norm{y_1y_2} \bol \mcox{\norm{ab}}.$
Consequently, we have
\begin{equation} \label{la ozenka snizu 3}
 \mcox{2(1 - \lapx{r})} \men \mcox{\norm{ab}} \men \norm{y_1y_2}.
\end{equation}

Using the intercept theorem, we obtain
\begin{equation}\label{la ozenka snizu 4}
  \norm{y_1y_2} = \frac{\norm{y_1y_2}}{\norm{oy_2}} = \frac{\norm{x_1x_2}}{\norm{x_0x_2}} =
  \frac{\norm{x_0x_2} - \norm{x_0x_1}}{\norm{x_0x_2}}  = 1 - \frac{r}{\norm{x_0x_2}} \men
  1 - \frac{r}{\xi_X}.
\end{equation}

By  inequalities \reff{la ozenka snizu 3} and \reff{la ozenka snizu 4}, we have 
$$ \mcox{2(1 - \lapx{r})} \men 1 - \frac{r}{\xi_X}$$
\end{prf}

It is easy to check that in a Hilbert space $H$ the following equality holds
$$
	\mcoi{H}{\tau} = 2 \sqrt{1 - (1 - \tau)^2}.
$$
Substituting this in  inequality  \reff{la ozenka snizu} 
and since  $\xi_H = 1,$  we obtain

$$
	\mco{H}{2r} = 1 - \frac{1}{2} \mcoi{H}{1 - {r}} \men \la^+_H (r).
$$

According to \reff{la_hilbert}, we have  that if $X$ is a Hilbert space,
 then the right hand  estimate  in  inequality   
 \reff{la ozenka snizu} is reached.

%

\end{document}